# The Wild Number Problem:
# math or fiction?

Philibert Schogt

> "Let's have a look then, shall we?" Dimitri said, as if he was a doctor about to examine a patient. I handed him the proof. He removed the paper clip and spread the pages out on his desk.
> " 'Theorem: the set of wild numbers is infinite'," he read the first line out loud. "My god, Isaac, why didn't you tell me you were working on this?"[1]

## 1 Escape from reality

In the above scene from my novel *The Wild Numbers*, the 35-year-old mathematician Isaac Swift believes he has found a solution to the famous Wild Number Problem and is showing it to his older colleague, the Russian mathematician Dimitri Arkanov. But what exactly is the Wild Number Problem, and why did I choose it as the theme for my novel, rather than some other mathematical problem?

Actually, my original plan was quite different.
 It was 1992, and I had just moved into a new apartment in Amsterdam. A friend of mine, who is a mathematician, was helping me install the electricity. While we were working, I told him that I had finally managed to construct a good plot for what was to become my first novel. I was sure he would appreciate its theme and was looking forward to his reaction.
 The main character was a troubled mathematician struggling with his mediocrity. For the past few years, his research has been heading absolutely nowhere. But then one day: eureka! In a sudden flash of insight he believes to have found the proof of Fermat's Last Theorem.
 "Not a good idea," my friend remarked.
 I was stunned. What I had just told him was only meant as an introduction. I hadn't even come to the main story-line or the most important characters. How could he be so quick to pass judgement?

---
[1] Ph. Schogt, *The Wild Numbers*, p. 101

What he objected to was the improbability of a mediocre mathematician suddenly solving what was generally considered to be the most famous unsolved problem in mathematics.

"He *thinks* he has solved the problem," I corrected him.

This nuance failed to impress him. The mathematics needed to tackle Fermat's Last Theorem was so complex, he explained, that only the very best minds would ever be able to find its proof. Unless my hero had miraculously turned into a genius overnight, he would have to be a complete idiot to think he had the answer. This was not the sort of protagonist my friend could identify with. Or sympathise with.

"But the average reader might be able to," I protested.

"Perhaps. But the average mathematician will never take someone like that seriously, not even for a moment."

His commentary was devastating. My aim was to write a novel about a mathematician that "the average reader" could understand and enjoy, but I wanted the story to be sufficiently true to life to appeal to professional mathematicians as well.

It had taken me months to develop my plan. Now, the very first mathematician who heard about it coldly dismissed it in less than half a minute. I was furious with him but was not in a position to let it show. After all, he was helping me with my apartment. We worked on in silence. I sat on the floor trying to put together a socket. Cramming too many wires into too small a box seemed just the right punishment for the grand gestures with which I had presented my plan. Meanwhile, with a grimace that struck me as somewhat Satanic, my friend began drilling holes in the wall. The shrill sound of the drill and the fine red powder that came spewing out of the holes were painful reminders of the damage he had just done to my ego.

When we were done, I did my best to thank him for his help. The lights were working, but the idea for my novel lay in ruins.

But the next day, I was blessed with one of those flashes of insight that I had wanted to bestow on my main character. The solution to my problems was so obvious that I wondered why I hadn't come up with it earlier. Instead of writing a novel about a real mathematical problem such as Fermat's Last Theorem, as a writer of fiction I was perfectly free to make up a mathematical problem of my own. And so I created the French mathematician Anatole Millechamps de Beauregard, who in his turn invented the wild numbers in 1823.

As I developed my new idea further, I realised that there were other drawbacks to the original plan which could now be avoided.

I am extremely lazy when it comes to doing research. If I had written my novel about Fermat, I would have had to go to the mathematics library, talk to various number theorists, strain my mathematical capabilities to the utmost in order to understand the latest developments.

Moreover, I would have had to consider introducing true historical figures into the story (more research!), keeping them in line with biographical facts and thus limiting my freedom to shape their personalities.

And then there was the problem that writers of popular science always run into: how do you present a complicated subject to the general public without oversimplifying or even distorting the facts?

Now that I had entered the realm of pure fiction, I no longer had to worry about such matters. Thanks to Beauregard and the wild numbers, the time I spent on research was reduced to a blissful near-zero.

**2 Fiction imitating math**

This does not mean to say that the mathematics in my novel could just be any old nonsense. Quite the contrary. One of the greatest challenges that I was now faced with was to make the Wild Number Problem seem as real as possible, appealing to the imagination of the general public and professional mathematicians alike.

The first thing that I needed to do was to have my French mathematician Anatole Millechamps de Beauregard define the problem, that is to say, I had to create the illusion that something was being defined. So I let him play a gambling game with his friends. They each deposited a sum of money, he would pose a mathematical riddle, and the first person to solve it would win the jackpot. One of these riddles came to be known as the Wild Number Problem:

> Beauregard had defined a number of deceptively simple operations, which, when applied to a whole number, at first resulted in fractions. But if the same steps were repeated often enough, the eventual outcome was once again a whole number. Or, as Beauregard cheerfully observed: "In all numbers lurks a wild number, guaranteed to emerge when you provoke them long enough". 0 yielded the wild number 11, 1 brought forth 67, 2 itself, 3 suddenly manifested itself as 4769, 4, surprisingly, brought forth 67 again. Beauregard himself had found fifty different wild numbers. The money prize was now awarded to whoever found a new one.[2]

The trick here, of course, was that I didn't specify what Beauregard's "deceptively simple operations" were.

The next step was to provide the Wild Number Problem with a history. It had to be old enough to count as a famous unsolved problem, but at the same time, to maintain the illusion that we were dealing with real mathematics, it was better to avoid going into too much historical detail. 1823 seemed just about the right starting point, allowing for three, maybe four important developments:

---

[2] *The Wild Numbers,* p. 34

1823:             Anatole Millechamps de Beauregard poses the Wild Number Problem in its original form.

1830s:            The problem is generalised: how many wild numbers are there? Do the same ones keep popping up, or are there infinitely many?

1907:             Heinrich Riedel ends speculations that perhaps all numbers are wild by proving that the number 3 is not. Later he extends his proof to show that there are infinitely many of such non-wild, or "tame" numbers.

early 1960s:      Dimitri Arkanov sparks renewed interest in the almost forgotten problem by discovering a fundamental relationship between wild numbers and prime numbers.

the present:      Isaac Swift finds a solution.

  When my main character Isaac Swift explains what the wild numbers are and recounts their history, he addresses the reader more or less directly. Here, I was able to imitate the tone of popular science books, intended to reassure readers that they can grasp the general idea without having to bother with the complexities of the subject. Only in this case, the existence of such complexities was entirely illusory.
  But elsewhere in my novel, when Isaac is working on the problem or discussing it with colleagues, the tone had to resemble the discourse of professional mathematicians, too technical for a layperson to understand. Here, I had to introduce some make-believe jargon into the narrative, with which I could pretend to be taking steps in a mathematical proof. As with historical details, the key once again was not to overdo it. Too much pseudo-mathematical mumbo jumbo would scare off the average reader and have mathematicians shaking their heads. In the end, I managed to limit the number of nonsense terms to five. Apart from the wild numbers themselves and their counterparts, tame numbers, I introduced the term "calibrator set", a mathematical tool developed by Isaac's older colleague Dimitri Arkanov that was useful in establishing a set's "K-reducibility"[3], a deep number theoretic property. And Isaac himself constructs "pseudo-wild numbers" – a type of number with a somewhat weaker definition than wild numbers – hoping to use them as a stepping stone in his proof.

---

[3] I later found out that there really is such a thing as K-reducibility, albeit in a branch of mathematics far-removed from number theory. The term must have caught my eye at some point in the past, lingering in my mind without my being aware of it.

With the help of these terms, I could now let my main character argue with himself while frenziedly pacing back and forth in his apartment:

> *Assuming there is a set of pseudo-wild prime numbers $Q_p$ that is infinite and K-reducible, find a correspondence between the elements $q_p$ and $w_p$ – wild primes – such that for every pseudo-wild prime there exists at least on wild prime...*
> "What do you mean *assume $Q_p$* is infinite and K-reducible? You are only shifting the problem!"[4]
>
> ...
>
> A calibrator set. A calibrator set. If only I found a suitable calibrator set![5]

As work on my novel progressed, I remember having crazy conversations with my mathematician friend, discussing various steps in the imaginary proof of a non-existent problem, and considering whether they were realistic enough. It was a bit like that famous scene in *Blow-up*, the film by Michelangelo Antonioni, where two people play tennis without a ball.

I don't think I would ever have had as much fun if I had stuck to my original plan and written a novel about Fermat's Last Theorem. In retrospect, I am deeply indebted to my friend, not only for helping me out with the electricity that day, but more importantly, for setting my imagination free.

## 3  Back to reality

When I completed the manuscript of *The Wild Numbers* in 1994, I left it lying on a shelf for quite some time before gathering the courage to send it to a literary agent. Several months later, he came with his diagnosis: great story, needs a sub-plot. More "human interest" was what he was looking for, fearing that the mathematics in my novel would scare off too many readers. Another agent responded in the same vein: she would have no trouble at all getting my novel published, if it were three times as long.

As an unpublished writer, I was inclined to take these professional opinions extremely seriously. But when I tried to follow up on the agents' suggestions, changing the narrative from first person to third person to allow for more non-mathematical detail, emphasising the love element, and so on, I felt that if anything, the story was weakened, not strengthened. The whole idea of my project had been to show that mathematics and its practitioners were interesting enough in themselves, and didn't need help from outside in the form of sub-plots to captivate the reader. And perhaps the brevity of my novel served a deeper purpose as well. One of the distinguishing features of mathematics is its aim to be clear and concise, to express everything in the simplest possible terms. Perhaps the considerations mentioned above - inventing a fictional problem to avoid the complexities of a real problem, limiting the amount of historical detail, keeping technical jargon to a minimum - could all be seen as an attempt to reduce the drama to its simplest form, thus reflecting the spirit of mathematics in a way that a thicker novel, fluffed up with interesting but superfluous detail, could not.

---
[4] *The Wild Numbers*, p. 82
[5] Ibid, p. 85, cf. Shakespeare's *Richard III*: "A horse. A horse. My kingdom for a horse!"

To make a long story short, I decided that my novel was fine the way it was. Unfortunately, this meant that my manuscript went back to gathering dust on a shelf.

Meanwhile, out in the real world, something truly amazing took place. When my original plan was still intact, I did wonder once or twice what would happen if Fermat's Last Theorem were actually proved while I was writing a story on the same subject. I had decided not to worry. The latest news back in 1992 was that some or other mathematician had published a proof, which, like so many others in the past, had turned out to be wrong. Considering the 350 years of failed attempts and modest steps in the right direction, surely it would be too much of a coincidence if a solution were found during the two or three years that I needed to write my novel.

And yet, this is exactly what happened.

In 1995, the English mathematician Andrew Wiles published his definitive proof of Fermat's Last Theorem. It was over 100 pages long, the result of seven years of hard work. If I had stuck to my old plan, I could have thrown my novel straight into the garbage can. Now that I had written about a fictional problem, at least my story was safe.

It was more than safe, actually. Partly thanks to Andrew Wiles' spectacular achievement, mathematics was becoming a hot topic in the media, and all sorts of books, films and plays were appearing on the market. Evidently, people were much more open to mathematics than the publishing industry had given them credit for. It was in this favourable climate that a Dutch editor read my manuscript and was willing to take the gamble. After I had translated the story into Dutch, it was finally published in 1998, under the title *De wilde getallen*. The original English version appeared in 2000. Since then, it has been translated into various other languages, including German, Greek and Italian.

## 4 Math imitating fiction

Whether I have succeeded in creating a credible and enjoyable piece of mathematical fiction is up to my readers to decide. But I consider it an encouraging sign that so many people have asked me if the wild numbers really exist. Friends of mine, who took my book along on their vacation with another family, told me that they had debated this issue on the beach one day. Their teenage daughter settled the question, insisting that her math teacher had discussed the wild numbers in class.

One literary critic, on the other hand, dismissed my novel because the mathematics in it was so evidently nonsensical. For a brief moment, I was dumbstruck. I had spent a great deal of time and effort making the wild numbers seem as real as possible. Where had I given myself away? Interestingly, the critic draws a correct conclusion from the wrong premises. His main objection was that the sequence of wild numbers mentioned in my book (11, 67, 2, 4769, 67) was too erratic to be realistic. In his opinion, I make matters worse by having my French mathematician Beauregard define his wild numbers with a "series of deceptively simple operations."

"I would certainly like to see those 'deceptively simple' operations!" the critic scoffs in his review.

But deceptive simplicity leading to erratic results is by no means peculiar to the wild numbers. In fact, this is a typical feature of many existing problems in number theory, one that inspired me to write my novel in the first place.

And I was happy to discover that *The Wild Numbers* in its turn appealed to a great number of mathematicians for the very same reason. Though seemingly fictional, the Wild Number Problem sparked a lively debate, centring on the issue of whether the wild numbers could be generated after all with a series of simple operations as described in my book, i.e. by operations which, when applied to an integer, would at first result in fractions, but upon sufficient iteration would once again produce an integer. Attempts to generate the exact sequence of numbers mentioned in my book were unsuccessful. But contrary to the literary critic's intuition, various mathematicians did come up with beautiful and indeed deceptively simple ways to produce similarly erratic integer sequences.

Here is one example, suggested by the Dutch mathematician Floor van Lamoen:

For a rational number p/q let f(p/q) = p*q divided by the sum of digits of p and q; a(n) is obtained by iterating f, starting at n/1, until an integer is reached, or if no integer is ever reached then a(n) = 0.

For example, for n=2:

$$2/1 \to 2/3 \to 6/5 \to 30/11 \to 330/5 = 66$$

Here are the first 48 terms of the sequence:

0, 66, 66, 462, 180, 66, 31395, 714, 72, 9, 5, 15, 3, 36, 42, 39, 2, 9, 45, 462, 12, 12, 90, 3703207920, 1692600, 84, 234, 27, 3043425, 74613, 6,

7930296, 264, 4290, 510, 315, 315, 73302369360, 1155, 3, 8, 239872017, 6, 4386, 1989, 18, 17740866, 499954980

To my delight, one outcome of these discussions was that the wild numbers were accepted as an entry in "The On-line Encyclopedia of Integer Sequences"[6], a huge data base developed by the American mathematician Neil Sloane, offering information on every thinkable kind of integer sequence. As an extra feature, all the sequences in the encyclopedia have been set to music, so I was even able to *listen* to the wild numbers! The various efforts of mathematicians to create Beauregard-like sequences also made it into the encyclopedia, being dubbed "pseudo-wild numbers".[7]

But the story does not end there. In 2004, I received a phone call from Jeff Lagarias, an American mathematician who was in Holland to attend a math conference. He had read my novel and was eager to meet me, so we had dinner together in a Thai restaurant in the centre of town. As it turned out, Lagarias was one of the world's leading experts on the unsolved 3n+1 problem, otherwise known as the Collatz conjecture. I was unfamiliar with the problem, but fortunately, he had no trouble at all explaining the principle to me, not because I am so intelligent, but because it is easy enough for a ten-year old to understand:

> Starting with any integer n, if even, divide by two, if odd, multiply by 3 and add 1. Repeat this process indefinitely. The conjecture states that no matter which integer we start out with, ultimately, the sequence will reach 1.
> 
> For instance, starting with $n = 6$, we get the sequence:
> 
> 6, 3, 10, 5, 16, 8, 4, 2, 1.
> 
> For $n = 11$, the sequence is:
> 
> 11, 34, 17, 52, 26, 13, 40, 20, 10, 5, 16, 8, 4, 2, 1.
> 
> If we consider the sequence for $n = 27$ and see how the numbers keep rising and falling before finally dropping down to 1, it becomes clear why the conjecture has been nicknamed "the hailstone problem":
> 
> 27, 82, 41, 124, 62, 31, 94, 47, 142, 71, 214, 107, 322, 161, 484, 242, 121, 364, 182, 91, 274, 137, 412, 206, 103, 310, 155, 466, 233, 700, 350, 175, 526, 263, 790, 395, 1186, 593, 1780, 890, 445, 1336,

---

[6] See Neil Sloane's On-line Encyclopedia of Integer Sequences™ (www.OEIS.org), entry A058883  
[7] Ibid, entries A058971, A058972, A058973, A058977, A058988, A059175

> 668, 334, 167, 502, 251, 754, 377, 1132, 566, 283, 850, 425, 1276, 638, 319, 958, 479, 1438, 719, 2158, 1079, 3238, 1619, 4858, 2429, 7288, 3644, 1822, 911, 2734, 1367, 4102, 2051, 6154, 3077, 9232, 4616, 2308, 1154, 577, 1732, 866, 433, 1300, 650, 325, 976, 488, 244, 122, 61, 184, 92, 46, 23, 70, 35, 106, 53, 160, 80, 40, 20, 10, 5, 16, 8, 4, 2, 1

When we came to talk about my book, Lagarias told me that the wild numbers had reminded him of his own work in so many ways, that he wanted to incorporate them into an article that he was planning on writing. I was greatly honoured, although I wasn't quite sure what to expect.

In 2006, Lagarias published two articles on the 3n + 1 problem, one in the *American Mathematical Monthly*, the other, co-authored with David Applegate in the *Journal of Number Theory*, in which he introduces the term "wild number", citing my novel as a source of inspiration. How Lagarias defines these wild numbers goes beyond the scope of our present discussion. For our purposes, it suffices to compare the following two statements:

> **Theorem.** The set of wild numbers is infinite.
> - Isaac Swift, *The Wild Numbers*

> **Theorem 3.1.** The semigroup of wild integers contains infinitely many irreducible elements (i.e., there are infinitely many wild numbers).
> - Jeffrey C. Lagarias, *Wild and Wooley numbers*,
>   American Mathematical Monthly 113 (2006), 97–108.

From now on, when readers ask me if the wild numbers really exist, I will be able to tell them the good news.